\documentclass[10pt]{amsart} \usepackage{amsthm} \usepackage{amsmath}
\usepackage{mystyle} \usepackage{amssymb} \usepackage{amsfonts}
\usepackage{hyperref} \usepackage{xy} \xyoption{all}
\newcommand{\la}{{\mathbb \lambda}} 
\newcommand{\Q}{{\mathbb Q}} 
\newcommand{\Z}{{\mathbb Z}} \newcommand{\F}{{\mathbb F}}
 \renewcommand{\div}{\mathrm{div}}
\newcommand{\red}{\mathrm{red}} \def\mclambda{\lambda} 
 \def\affmod#1#2{{\text{($\text{mod}_#1$ $#2$)}}}

\numberwithin{equation}{section}

\begin{document}
\title{On Shafarevich-Tate groups and the Arithmetic of Fermat Curves}
\dedicatory{To Sir Peter Swinnerton-Dyer on his 75th birthday.}
\author{William G. McCallum} \author{Pavlos Tzermias}

\address{Department of Mathematics, P.O. Box 210089, 617 N. Santa
  Rita, The University of Arizona, Tucson, AZ 85721-0089}
\email{wmc@math.arizona.edu}

\address{Department of Mathematics, University of Tennessee,
  Knoxville, TN 37996-1300} \email{tzermias@math.utk.edu}

\maketitle

\section{Introduction}\label{sec:introduction}
Let $\Q$ denote the field of rational numbers and $\overline{\Q}$ a
fixed algebraic closure of $\Q$.  For a fixed prime $p$ such that $p
\geq 5$, choose a primitive $p$-th root of unity $\zeta$ in
$\overline{\Q}$ and let $K=\Q(\zeta)$.  If $a$, $b$ and $c$ are
integers such that $0 < a, \ b, \ a+b < p$ and $a+b+c=0$ let
$F_{a,b,c}$ denote a smooth projective model of the affine curve
\begin{equation}\label{eq:1}
y^{p} = x^{a} (1-x)^{b}
\end{equation}
and let $J_{a,b,c}$ be the Jacobian of $F_{a,b,c}$.  Then $J_{a,b,c}$
has complex multiplication induced by the birational automorphism
$(x,y) \mapsto (x, \zeta y)$ of $F_{a,b,c}$.  Let $\la$ denote the
endomorphism $\zeta - 1$ of $J_{a,b,c}$.  Note that $\la^{p-1}$ is, up
to a unit in $\Z[\zeta]$, multiplication by $p$ on $J_{a,b,c}$.

We are interested in the Shafarevich-Tate group of $J_{a,b,c}$ over
$K$, which we denote simply by $\Sh$.  In \cite{mccallum:1988}, the
first author studied the restriction of the Cassels-Tate pairing
\begin{equation}\label{eq:2}
\Sh[\la] \times \Sh[\la] \longrightarrow \Q/\Z
\end{equation}
and showed that $\Sh[\la]$ is non-trivial in certain cases depending
on the reduction type of the minimal regular model of $F_{a,b,c}$ over
$\Z_p[\zeta]$.  The purpose of this paper is to extend those results
by carrying out higher descents, and to derive some consequences for
the arithmetic of Fermat curves using the techniques of the second
author.

\begin{figure}\label{fig:reductiontype}
\caption{Reduction types of $F_{a,b,c}$}
\begin{picture}(400,140)(10,0)
\put(85,55){\footnotesize multiplicity 2}
\put(90,110){Tame }
\put(250,110){Wild }
\thicklines
\put(10,75){\parbox{.6in}{\footnotesize proper transform}}
\put(30,50){\line(1,0){150}}
\thinlines
\put(50,40){\line(0,1){50}}
\put(70,40){\line(0,1){50}}
\put(100,75){\circle*{2}}
\put(115,75){\circle*{2}}
\put(130,75){\circle*{2}}
\put(160,40){\line(0,1){50}}
\put(50,35){\centering{$\underbrace{\phantom{aaaaaaaaaaaaaaaaaaaaa}}$}}
\put(95,20){\footnotesize $p+1$}
\qbezier{(250,100)(300,55)(250,10)}
\qbezier{(300,100)(250,55)(300,10)}
\put(227,55){\line(1,0){100}}
\put(227, 42){\parbox{.6in}{\footnotesize proper transform}}
\end{picture}
\end{figure}
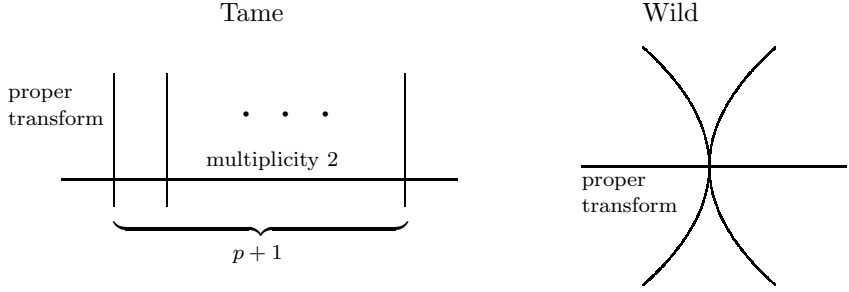

First we recall the main result of \cite{mccallum:1988}.  The possible
reduction types for $F_{a,b,c}$ are as shown in
Figure~\ref{fig:reductiontype} \cite{mccallum:1982}, with the proper
transform of the special fiber of the model~(\ref{eq:1}) indicated.
The wild type is further divided into split and non-split, according
to whether the two tangent components are defined over the finite
field $\FF_p$ or conjugate over a quadratic extension.  The reduction
type can be computed as follows.  For a rational number $x$ of
$p$-adic valuation $0$, let $q(x)=(x^{p-1}-1)/p$, viewed as an element
of $\F_p$. Then $F_{a,b,c}$ is
$$
\begin{array}[c]{ll}
\text{ tame }& \text{if $-2 a b c q(a^a b^b c^c) = 0$}\\
\text{ wild split} & \text{if 
$-2 a b c q (a^a b^b c^c) \in \F_p^{\times 2}$} \\
\text{ wild non-split} & \text{if  
$-2 a b c q(a^a b^b c^c) \notin \F_p^{\times 2}$}
\end{array}
$$
Let $M_K$ be the set of finite places of $K$ and let $w$ denote the
unique place of $K$ above $p$. Define
\begin{equation}
U = \{ x \in K^{\times}/K^{\times p} : \mbox{$v(x) \equiv 0
  \pmod{p}$ for all $v \in M_K$} \},\qquad V =
K_w^{\times}/K_w^{\times p}.\label{eq:36}
\end{equation}
Let $\pi$ be the uniformizer of $K_w$
defined by 
\begin{equation}
\text{$\pi^{p-1}=-p$ and $\frac{\pi}{1-\zeta} \equiv 1
\pmod{w}$.}\label{eq:35}
\end{equation}
If $\kappa : \Gal(K/\Q) \to \Z_p^{\times}$ is the Teichm\"{u}ller
character, let $V(i)$ denote the intersection of the $\kappa^{i}$-th
eigenspace of $V$ with the subgroup of $V$ generated by units
congruent to $1$ modulo $\pi^i$ (thus $V(i)$ is one-dimensional if $2
\le i \le p$).
\begin{theorem}[\cite{mccallum:1988}]\label{thm:old}
  Suppose that $F_{a,b,c}$ is wild split, $p \equiv 1 \pmod{4}$, and
  the image of $U$ is non-trivial in both $V((p-1)/2)$ and
  $V((p+3)/2)$.  Then $\Sh[\lambda]/\lambda \Sh[\lambda^2] \iso
  \ZZ/p\ZZ\oplus \ZZ/p\ZZ$.
\end{theorem}

The condition on $U$ is satisfied if $p\nmid B_{(p-1)/2}B_{(p+3)/2}$,
where $B_k$ is the $k$-th Bernoulli number.  As noted in
\cite{mccallum:1988}, the technique used to prove
Theorem~\ref{thm:old} applies to the pairing
\begin{equation}\label{eq:3}
\Sh[\la^2] \times \Sh [\la] \longrightarrow \Q/\Z
\end{equation}
and yields information about $\Sh[\la^2]$. 
\begin{theorem}\label{thm:free}
  Suppose that either of the following conditions is satisfied:
  \renewcommand{\theenumi}{\alph{enumi}}
\begin{enumerate}
\item $F_{a,b,c}$ is wild-split and $p \equiv 3
  \pmod{4}$\label{item:1}
\item $F_{a,b,c}$ is wild non-split or tame and the image of $U$ in
  either $V((p+1)/2)$ or $V((p+3)/2)$ is trivial. \label{item:2}
\end{enumerate}
Then $\Sh[\lambda^2]/\lambda \Sh[\lambda^3] = 0$, that is,
$\Sh[\la^3]$ is a free module over $\Z[\zeta]/(\la^3)$.
\end{theorem}

As discussed in \cite{mccallum:1988}, the hypothesis on $U$ in
condition~(\ref{item:2}) of the theorem is quite mild, since for $U$
to be nontrivial in $V(k)$ with $k > 1$ and odd requires that $p$
divides $B_{p-k}$.

\begin{corollary}\label{thm:corollary}
  If one of conditions~(\ref{item:1}) or (\ref{item:2}) of
  Theorem~\ref{thm:free} is satisfied, and if $|\Sh[p^{\infty}]| <
  p^{3}$, then $\Sh[p^{\infty}] = 0$.
\end{corollary}

It is natural to ask which occurs more often under the conditions of
Theorem~\ref{thm:free}, $|\Sh[p]| = 0$ or $|\Sh[p]| \ge p^3$. To
explore this question, we compute
\begin{equation}\label{eq:4}
\Sh[\la^3] \times \Sh [\la] \longrightarrow \Q/\Z.
\end{equation}
\begin{theorem}\label{thm:nontrivial}
  Suppose that $p\ge19$ is regular, $p \equiv 3 \pmod{4}$, $F_{a,b,c}$
  is tame or wild non-split and
  \begin{equation*}
    q(a^ab^bc^c)^3+ a b c B_{p-3} \not\equiv 0 \pmod{p}.
  \end{equation*}
  Then $\Sh[\la^3] \neq 0$ (and hence, by
  Corollary~\ref{thm:corollary}, $|\Sh[p^\infty]| \ge p^3$).
\end{theorem}

The next proposition shows that, in certain cases, one can combine
Theorems~\ref{thm:free} and~\ref{thm:nontrivial} to describe the exact
structure of $\Sh[p^{\infty}]$:

\begin{theorem}\label{prop:exact-p-primary} 
  Suppose that $p$, $a$, $b$ and $c$ are chosen so that the hypotheses
  of both Theorems~\ref{thm:free} and~\ref{thm:nontrivial} are
  satisfied. If, in addition, the free $\Z[\zeta]/(\la^3)$-module
  $\Sh[\la^3]$ has rank 2, then
  $$\Sh[p^{\infty}]=\Sh[\la^3] \iso (\Z[\zeta]/(\la^3))^2. $$
\end{theorem}

\noindent In Section~\ref{sec:tame-reduction} we establish the following
application of the above results:

\begin{theorem}\label{cor:diophant}
  Let $p=19$, $a=7$, $b=1$. Then
\begin{enumerate}
\item $\Sh[p^{\infty}] \iso (\Z[\zeta]/(\la^3))^2$. \label{diophant:1}
\item The Mordell-Weil rank of $J_{7,1,-8}$ over $\Q$ equals 1.
           \label{diophant:2}
\item The only quadratic points (i.e.\ algebraic points whose field of
  definition is a quadratic extension of $\Q$) on the Hurwitz-Klein
  curve $F_{7,1,-8}$ and also on the Fermat curve
  $X^{19}+Y^{19}+Z^{19}=0$ are the ones described by Gross and
  Rohrlich in \cite{gross-rohrlich:1978}.
           \label{diophant:3}
\end{enumerate} 
\end{theorem}

\noindent 
We also note that, by combining Theorem~\ref{thm:nontrivial} with
Faddeev's bounds in \cite{faddeev:1961}, one gets that the
Mordell-Weil rank (over $\Q$) of any tame or wild non-split quotient
of the Fermat curve $F_{19}$ or $F_{23}$ is at most 2.

Lim \cite{lim:1995} has also stated a result attempting to improve on
\cite{mccallum:1988} in certain cases. However, in
Section~\ref{sec:tame-reduction}, we show that the hypotheses of
Propositions A and B of \cite{lim:1995} are never simultaneously
satisfied.

\section{Formulas for the pairings}\label{sec:formulas-pairings}
We recall the situation and notation of\cite{mccallum:1988}. For $\phi
\in {\mathcal O}_K$ and $F$ a field containing $K$, we denote by
$\delta = \delta_{\phi,F}$ the coboundary map $J(F) \to H^1(F,
J[\phi])$.  We have the $\phi$-Selmer group $S_\phi \subset H^1(K,
J[\phi])$, defined to be the subgroup whose specialization to each
completion $K_v$ of $K$ lies in the image of $\delta_{\phi,K_v}$, and
which sits in an exact sequence
\begin{equation*}
  0 \to J(K)/\phi J(K) \to S_\phi \to \Sh[\phi] \to 0.
\end{equation*}
For $\phi, \psi \in \End(J)$, we have a pairing
 \begin{equation}
 \label{eq:SelmerPairing}
S_\phi \times S_{\hat\psi} \rightarrow
\mathbb{Q}/\mathbb{Z},
\end{equation}
described in \cite{mccallum:1988}, which is a lift of the restriction
of the Cassels pairing to $\Sh[\phi]\times\Sh[\hat\psi]$.  An
expression for the pairing~\eqref{eq:SelmerPairing} is given in
\cite{mccallum:1988}, under a certain splitting hypothesis.

We use this formula to derive formulas for the pairings~(\ref{eq:3})
and~(\ref{eq:4}). The formula for~(\ref{eq:3}) is a straightforward
consequence of Theorem 2.6 in \cite{mccallum:1988}; the formula
for~(\ref{eq:4}) takes more work. The point is that $J[\lambda^3]
\subset J(K)$ \cite{greenberg:1981}, so that it is possible to express
the pairings (\ref{eq:2}) and (\ref{eq:3}) as purely local pairings at
$w$, as explained in \cite{mccallum:1988}. However, by
\cite{greenberg:1981} and \cite{kurihara:1992}, the
$\lambda^4$-torsion on $J_{a,b,c}$ is not in general defined over $K$,
introducing an essentially global aspect to the calculation
of~\eqref{eq:4}.

For technical reasons, it is convenient to replace $\lambda$ with an
endomorphism (which we also denote by $\lambda$) that is congruent
modulo $\lambda^5$ to the uniformizer $\pi$ defined by~\eqref{eq:35}, 
since then
\begin{equation*}
  \mclambda^\delta \equiv \mclambda^{\kappa(\delta)} 
  \pmod{\mclambda^5},\quad \delta \in
  \Gal(K/\mathbb{Q}). 
\end{equation*}
In particular, we have $\hat\mclambda \equiv -\mclambda$ modulo
$\mclambda^5$, and we will often replace $\hat\mclambda$ with
$-\mclambda$ in what follows without mention, in cases where we are
dealing with a module killed by $\mclambda^5$. Furthermore, it
suffices to prove Theorems~\ref{thm:free} and~\ref{thm:nontrivial}
with this new choice of $\mclambda$.  Since $\mclambda/\hat\mclambda$
is a unit, $\Sh[\mclambda] = \Sh[\hat\mclambda]$, and we can proceed
by computing the pairing \eqref{eq:SelmerPairing} with $\phi =
\mclambda^k$ and $\psi = \hat\mclambda$.

The local formula for the Cassels-Tate pairing is expressed in terms
of certain local descent maps as follows.  Given a $p$-torsion point
$Q$ in $J(\overline K)$ we denote by $D_Q$ a divisor defined over
$K(Q)$ representing $Q$ and by $f_Q$ a function on $F_{a,b,c}$ whose
divisor is $pD_Q$. Evaluating $f_Q$ on divisors induces a map
$\iota_Q:J(F) \to F^\times/F^{\times p}$ for any field $F$ containing
$K(Q)$.

By \cite{greenberg:1981}, $K(J[\mclambda^3]) = K$ and
\begin{equation}
K(J[\mclambda^4]) = L
= K(\eta_{p-3}^{1/p}), \label{eq:5}
\end{equation}
where $\eta_{p-3}$ is a generator for the $\kappa^{p-3}$-eigenspace of
the cyclotomic units in $K$. Let $\tilde \Delta \subset
\Gal(L/\mathbb{Q})$ be a subgroup projecting isomorphically to
$\Gal(K/\mathbb{Q})$.  For $i = 1, 2, 3, 4$, we choose points $P_{i}$
of order $\mclambda^{i}$ on $J$ and a generator $\sigma$ for
$G=\Gal(L/K)$ such that
\begin{enumerate}
\item $P_1$ is the point represented by the divisor $(0,0)-\infty$
\item $\mclambda P_i = P_{i-1}$, $i = 2, 3, 4$
\item $P_i$ is an eigenvector for the action of $\tilde \Delta$ with
  character $\kappa^{1-i}$
\item $\sigma P_4 = P_4 + P_1$.
\end{enumerate}
For $i \leq 4$, let $e_{\lambda^i}(P,Q)$ be the $\lambda^i$ Weil pairing on
$J[\lambda^i]$. We have an 
isomorphism $J[\mclambda^i] \iso \mu_p^i$ defined over
$K(P_i)$ (and thus over $K$ for $i \leq 3$), give by 
\begin{equation}
  \label{eq:37}
  Q \mapsto (e_{\lambda^i}(Q, P_1), \dots , e_{\lambda^i}(Q, P_i)).
\end{equation}
With this identification we have, by \cite[Lemma
2.2]{mccallum:1988},
\begin{equation*}
  \delta_i = \delta_{\mclambda^i, K(P_i)} = 
   \iota_{P_1} \times \cdots \times \iota_{P_i}.
\end{equation*}
Since
$J$ has good reduction outside $p$ and $\mclambda$ has degree $p$, we
can regard $S_{\mclambda^i}$ as a subgroup of $H^1(K(p)/K,
J[\mclambda^i])$, where $K(p)/K$ is the maximal extension of $K$
unramified outside $p$. Then, for $i\leq 3$,  $H^1(K(p)/K, J[\mclambda^i])$ is
identified with the subgroup $U^i$ of $(K^\times/K^{\times p})^i$,
where $U$ is as defined in~\eqref{eq:36}.
For $a, b \in K_w^\times$, denote by $(a,b)$ the Hilbert symbol.
\begin{proposition}\label{thm:pairingformula}
  Let $a \in S_{\mclambda^2}$, $b \in S_{\hat\mclambda}$, $a_w =
  \delta(x_w)$, $x_w \in J(K_w)$. Then
  $$\zeta^{p\langle a, b \rangle_2}
  =(\iota_{P_3}(x_w), b_w).$$
\end{proposition}

\begin{proof}
  This follows from Theorem~2.6 in \cite{mccallum:1988}, with $\phi =
  \mclambda^2$ and $\psi = \mclambda$.
\end{proof}
For a number field $F$ we denote by $\mathcal{O}_F^{\prime}$ the ring
of $p$-integers in $F$. If $F \subset K(p)$, then we have exact
sequences
\begin{equation}
  \label{eq:6}
  0 \rightarrow \mathcal{O}_F^{\prime
  \times}/\mathcal{O}_F^{\prime \times p} \rightarrow H^1(K(p)/F,
  \mu_p) \rightarrow C[p] \rightarrow 0
\end{equation}
and
\begin{equation}\label{eq:7}
  0 \rightarrow C/pC \rightarrow H^2(K(p)/F, \mu_p) \rightarrow
  \Br(F)[p] \rightarrow 0,
\end{equation}
where $C$ is the ideal class group of $F$.

\begin{lemma}\label{lemma:lifting}
  Each element of $H^1(K(p)/K, J[\mclambda^k])$ lifts to $H^1(K,
  J[\mclambda^{k+1}])$. Furthermore, if $p$ is regular, it lifts to
  $H^1(K(p)/K, J[\mclambda^{k+1}])$.
\end{lemma}
\begin{proof}
  Let $a \in H^1(K(p)/K, J[\mclambda^k])$, and let $\delta a \in
  H^2(K(p)/K, J[\mclambda])$ be  the coboundary of $a$ for
  the sequence
\begin{equation}
0 \to J[\mclambda] \to J[\mclambda^{k+1}] \to J[\mclambda^k] \to
0.
\end{equation}
Then the inflation of $\delta a$ in $H^2(K, J[\mclambda]) \iso H^2(K,
\mu_p) = \Br(K)[p]$ has zero invariant at every place not dividing
$p$. Thus it is zero, by the Brauer-Hasse-Noether theorem (since there
is only one place of $K$ dividing $p$). For the second statement , we
argue in the same way, using~\eqref{eq:7}.
\end{proof}

We recall the definition of $\langle , \rangle_3$. Let $a \in
S_{\mclambda^3}$ and $b \in S_{\hat \mclambda}$. Lift $a$ to an
element $a_1$ of $H^1(K, J[\mclambda^{4}])$ (possible by
Lemma~\ref{lemma:lifting}). For each place $v$ of $K$, lift $a_v$ to
an element $a_{v,1}$ that is in the image of $\delta$. Then $a_{1,v} -
a_{v,1}$ is the image of an element $c_v \in H^1(K_v, J[\mclambda])$,
and
\begin{equation*}
 \langle a, b \rangle = \sum_v c_v \cup b_v
\end{equation*}
where the cup product is with respect to the local pairing
$$H^1(K_v, J[\mclambda]) \times H^1(K_v, J[\hat\mclambda]) \to
\QQ/\ZZ.$$
If $p$ is regular, $L/K$ is totally ramified at $w$, and
there is a unique extension of $w$ to $L$, which we also denote by
$w$.  Our calculation of $\langle , \rangle_3$ uses the following
lemma.
\begin{lemma}\label{thm:lemma-let-l}\label{thm:lemma-suppose-that}
  Suppose $p>5$ is regular, and let $L$ be as in~\eqref{eq:5}. Then
  \begin{enumerate}
  \item  the
  map $H^1(K(p)/L, \mu_p) \rightarrow H^1(L_w, \mu_p)$ is injective
\item the norm map $N_{L/K}:
  \mathcal{O}_L^\times \rightarrow \mathcal{O}_K^\times$ is
  surjective
\item 
$H^1(K(p)/L, \mu_p)^G(i) = 0$ if $i$ is odd and $i \neq 1$, or if $i=p-1$.
  \end{enumerate}
\end{lemma}

\begin{proof}
  Let $H_K$ (resp.  $H_L$) be the Hilbert class field of $K$ (resp.
  $L$). Since $L/K$ is unramified outside $p$, and there is only prime
  of $L$ above $p$, it follows that $\Gal(H_L/L)/(\sigma-1) \iso
  \Gal(H_K/K)$. Therefore $p$ does not divide the order of
the class group   $C_L \iso \Gal(H_L/L)$. The
  injectivity statement follows, since anything in the kernel would
  generate an unramified Kummer extension of $L$ of degree $p$. 
Furthermore, every unit of
  $K$ is a local norm everywhere except possibly at the prime above
  $p$, and therefore is a local norm there also by the product
  formula. Thus it is a global norm. The surjectivity of the norm map
  follows by a standard argument using
  $\Gal(L/K)$ cohomology of the sequences
\begin{equation*}
1 \rightarrow \mathcal{O}_L^\times \rightarrow L^\times \rightarrow
  P_L \rightarrow 1
\end{equation*}
and
\begin{equation*}
  1 \rightarrow P_L \rightarrow I_L \rightarrow C_L \rightarrow 1,
\end{equation*}
where $I_L$ and $P_L$ are the groups of ideals and principal ideals
respectively.
 Finally, by~\eqref{eq:6},
  $H^1(K(p)/K, \mu_p) = \mathcal{O}_K^{\prime
    \times}/\mathcal{O}_K^{\prime \times p}$ and $H^1(K(p)/L, \mu_p) =
  \mathcal{O}_L^{\prime \times}/\mathcal{O}_L^{\prime \times p}$. 
Furthermore, the
  cokernel of $\mathcal{O}_K^\times/\mathcal{O}_K^{\times p}$ in
  $(\mathcal{O}_L^\times/\mathcal{O}_L^{\times p})^G$ is $H^2(L/K,
  \mu_p) \iso \ZZ/p\ZZ$, with $\Gal(K/\mathbb{Q})$ acting via
  $\kappa^{p-3}$, since it acts on $G$ via $\kappa^{3}$. Since $p-3$
  is even and $(\mathcal{O}_K^\times/\mathcal{O}_K^{\times p})(i) = 0$
  if $i$ is odd and $i\neq 1$, or if $i = p-1$, the third statement of the lemma
  follows.
\end{proof}
Let $N' = \sum_{i=1}^{p-1} i\sigma^i$.
\begin{proposition}\label{thm:prop-supp-p}
   Suppose $p$ is regular. Let $a \in S_{\mclambda^3}$, $b \in
   S_\mclambda$, $a_w = \delta(x_w)$, $x_w
   \in J(K_w)$.  Write $\mclambda^2_*a$, regarded as an element of
   $\mathcal{O}_K^\times/\mathcal{O}_K^{\times p}$, as $N_{L/K}\epsilon$ for
   some $\epsilon \in {\mathcal{O}}_L^\times$.  Then there exists a
   $\mclambda^4$-torsion point $P_4$, and an element $c_w \in
   K_w^\times$ such that
   $$
   \zeta^{p \langle a, b \rangle_3} = (c_w, b_w)
   $$
   and the image of $c_w$ in $L_w$ satisfies
   $$
   c_w = \iota_{P_4}(x_w)^{-1}N'\epsilon.
   $$
\end{proposition}

\begin{proof}
  Since $p$ is regular, $S_\lambda \subset
  \mathcal{O}_K^\times/\mathcal{O}_K^{\times p}$, so $b$ has zero
  component in the $\kappa^i$ eigenspace if $i$ is odd and if $i=p-1$ (the possibility $i=1$ is eliminated
  by local conditions, see \cite[6.2]{mccallum:1988}). Thus, by
  equivariance properties of the Hilbert pairing, we may replace $c_w$
  by its projection onto the $\sum_i (K_w^\times/K_w^{\times p})(i)$,
  with the sum over odd $i > 1$.  Consider the sequence
\begin{equation} 
\label{eq:8}
0 \to J[\mclambda] \to J[\mclambda^4] \to J[\mclambda^3] \to 0, 
\end{equation}
and the commutative diagram with exact rows
  \begin{equation*}
\xymatrix{0 \ar[r] & H^1(K(p)/L, J[\mclambda])^G \ar[r] & H^1(K(p)/L, J[\mclambda^4])^G
  \ar[r]^{\mclambda_*} & H^1(K(p)/L, J[\mclambda^3])^G \\ 
 & H^1(K(p)/K, J[\mclambda]) \ar[r] \ar[u]^{\res_{L/K}} & 
H^1(K(p)/K, J[\mclambda^4])
  \ar[u]^{\res_{L/K}}\ar[r]^{\mclambda_*} & H^1(K(p)/K, 
J[\mclambda^3])\ar[u]^{\res_{L/K}} }
  \end{equation*}
  The top row is exact because (\ref{eq:8}) splits over $L$, and hence
  the sequence
  $$
  0 \to H^1(K(p)/L, J[\mclambda]) \to H^1(K(p)/L, J[\mclambda^4])
  \to H^1(K(p)/L, J[\mclambda^3])
  $$
  is exact.  By Lemma~\ref{lemma:lifting} $a$ lifts to an element
  $a_1 \in H^1(K(p)/K, J[\mclambda^4])$. Let $a'_1 \in H^1(K(p)/L,
  J[\mclambda^4])^G$ be any lift of $\res_{L/K} a$ ($\res_{L/K} a_1$
  itself is one such).  Then
\begin{equation}\label{eq:9}
\res_{L/K} a_1 = a'_1 \eta,\quad \eta \in H^1(K(p)/L, \mu_p)^G.
\end{equation}
We now construct a candidate for $a'_1$. Under the
identification~\eqref{eq:37} between $J[\lambda^i]$ and $\mu_p^i$, the
map $\lambda^{i-1}: J[\lambda^i] \rightarrow J[\lambda]$ corresponds to
projection on the first component. Hence, in the  identification
$H^1(K, J[\mclambda^3]) = (K^\times/K^{\times p})^3$,
$a$ corresponds to an element $(x_1, x_2, x_3) \in (K^\times/K^{\times
  p})^3$, and $\mclambda^2_*a = x_1$.  Furthermore, in the identification
$$
H^1(L, J[\mclambda^4]) \iso (L^\times/L^{\times p})^4,
$$
the action of $\sigma$ on $H^1(L, J[\mclambda^4])$ is
intertwined with
$$
(t_1, t_2, t_3, t_4) \mapsto (t_1^\sigma , t_2^\sigma,
t_3^\sigma, t_4^\sigma t_1^\sigma), \quad t_i \in L^\times/L^{\times p}.
$$
Thus $(t_i)$ is fixed by $G$ if
$$
t_i^\sigma = t_i, \quad i = 1, 2, 3, \quad \text{and} \quad
t_4^{\sigma-1} = t_1^{-1}.
$$
By hypothesis, $x_1 = N_{L/K}\epsilon$, $\epsilon \in {\mathcal
  O}_L^\times$. Then
\begin{equation}\label{eq:10}
a'_1 = (x_1,
  x_2, x_3, N'\epsilon)
\end{equation}
is an equivariant lift of $(x_1, x_2, x_3)$.
  
Now let $a_{w,1}$ be the local lift of $a$ given by $a_{w,1} =
\delta_4(x_w)$.  Then
\begin{equation}\label{eq:11}
\res_{L_w/K_w} a_{w,1} = (x_1,  x_2, x_3, \iota_{P_4}(x_w)).
\end{equation}
Thus, from equations~\eqref{eq:9}, \eqref{eq:10}, and \eqref{eq:11},
we get
\begin{equation*}
\res_{L_w/K_w}(c_w) = \res_{L_w/K_w}(a_{1,w} - a_{w,1}) =
  \iota_{P_4}(x_w)^{-1}\eta N' \epsilon.
\end{equation*}
Since $\eta \in H^1(K(p)/L, \mu_p)^G$, its projection onto an
eigenspace $(K_w^\times/K_w^{\times p})(i)$ with $i > 1$ odd is
trivial. Thus, by the remarks at the beginning of the proof, we may
ignore the contribution of $\eta$, and the proposition follows.
\end{proof}

 \section{The local approximation}\label{sec:doing-local-appr}
 Let $P_i$ be as in the previous section, $i = 1, 2, 3, 4$, let $D_i$
 be a divisor on $F_{a,b,c}$ representing $P_i$, and let $f_i$ be a
 function whose divisor is $pD_i$. Take $D_i$ and $f_i$ to be defined
 over $K = \mathbb{Q}(\zeta)$ if $i = 1, 2, 3$ and over $L =
 K(\eta_{p-3}^{1/p})$ if $i = 4$.  The maps $\iota_{P_i}$ in
 Propositions~\ref{thm:pairingformula} and~\ref{thm:prop-supp-p} are
 computed by evaluating $f_i$ on certain divisors.  We use the
 approximation technique in \cite{mccallum:1988} to find expansions
 for $f_i$ on $p$-adic discs in $F_{a,b,c}$. Given a function $f$
 whose divisor is divisible by $p$ we approximate $f$ on an affinoid
 $Y$ in $F_{a,b,c}$ using the fact that
\begin{equation}
  \label{eq:12}
  \frac{df}{f} \equiv \omega \quad \affmod{Y}{p}.
\end{equation}
for some holomorphic differential $\omega$ on $F_{a,b,c}$
(\cite{mccallum:1988}, Theorem 5.2).  For general facts about rigid
analysis, we refer the reader to \cite{bosch-guentzer-remmert:1984}.

We recall the notion of congruence used in~\eqref{eq:12}.  If $Y$ is a
one-dimensional affinoid defined over an extension $F$ of $\QQ_p$ with
uniformizer $\pi_F$, we let $A(Y)$ be the ring of rigid analytic
functions on $Y$, $M(Y)$ the quotient field of $A(Y)$, and $D(Y)$ the
module of Kahler-differentials of $M(Y)$. We define
sub-$\mathcal{O}_F$-modules
\begin{eqnarray*}
  A^0(Y) &=& \{ f \in A(Y): \text{$f(x) \leq  1$ for all $x \in
  Y(\mathbb{C}_p)$}\}\\
  M^0(Y) &=& \{ f/g: f \in A^0(Y), g \in A^0(Y) \setminus \pi_F
  A^0(Y)\}\\
  D^0(Y) &=& \{ f\,dg: f, g \in M^0(Y)\}.
\end{eqnarray*}
If $f, g \in A(Y)$, $c \in F$, we say $f \equiv g \affmod{Y}{c}$ if
$(f-g) \in c A^0(Y)$, and similarly we define the notion of congruence
 on $Y$ in $M(Y)$ and $D(Y)$.  In order to deduce
from~\eqref{eq:12} information about power series expansions of $f$ on
closed discs in $Y$, we need the following lemmas.

\begin{lemma}\label{thm:lemma-if-z}
  Suppose that $Y$ is a one-dimensional affinoid over a finite
  extension $F$ of $\mathbb{Q}_p$, $Y$ has good reduction, and $Z$ is
  an affinoid contained in $Y$, isomorphic to a closed disc. If
  $\omega \in D^0(Y)$ is a differential with at worst simple poles on
  $Y$ that is regular on $Z$, then $\omega \in D^0(Z)$.
\end{lemma}

\begin{proof}
  Since $Z$ is isomorphic to a closed disc, it is contained in a
  residue class $U$ of $Y$ (or is equal to $Y$, in which case there is
  nothing to prove). It is clear from the definitions that $D^0(Y)|_U
  = D^0(U)$, hence $\omega \in D^0(U)$. Furthermore, since $Y$ has
  good reduction, $U$ is isomorphic to an open disc. Choose a
  parameter $t$ for $U$ such that $Z$ is the disc $|t|\leq |c| < 1$
  for some $c \in F $, and write
  \begin{equation*}
    \omega = g\,dt + \sum_{i=1}^n \frac{a_i}{t-b_i}\,dt,\quad g \in
    \mathcal{O}_F[[t]], a_i, b_i \in \mathcal{O}_F, |c| < |b_i| < 1.
  \end{equation*}
  Expanding the polar terms in powers of $t/b_i$ and setting $t = cs$,
  we get $\omega = f \, ds$ for some $f \in \mathcal{O}_F[[s]]$. Since
  $s$ is a parameter on $Z$, this proves the lemma.
\end{proof}

\begin{lemma}\label{thm:YtoX}
  Suppose that $f$ is a function whose divisor is divisible by $p$.
  Let $Y$ be an affinoid with good reduction contained in $F_{a,b,c}$
  and let $Z$ be a $p$-adic disc contained in $Y$ such that there is a
  function on $F_{a,b,c}$ restricting to a parameter on $Z$.  If
  $\omega$ satisfies the congruence (\ref{eq:12}) $\affmod{Y}{p}$,
  then it satisfies the same congruence $\affmod{Z}{p}$.
\end{lemma}

\begin{proof}
  With notation as in \cite{mccallum:1988}, we have
  $$
  \frac{df}{f} = \omega + p \eta, \quad \eta \in D^0(Y).
  $$
  Let $g$ be a function on $F_{a,b,c}$ such that $f/g^p$ is regular
  on $Z$ (we can construct $g$ using a parameter on $Z$ as in the
  hypotheses).  Since a suitable scalar multiple of $g$ is in
  $M^0(Y)$, $\dlog g \in D^0(Y)$. Thus $\eta - \dlog g \in D^0(Y)$ and
  is also regular on $Z$, and hence is in $D^0(Z)$ by
  Lemma~\ref{thm:lemma-if-z}. Thus
  $$
  \frac{df}{f} \equiv \frac{df}{f} - p\frac{dg}{g} = \omega +
  p(\eta - \frac{dg}{g}) \equiv \omega \quad \affmod{Z}{p}.
  $$
\end{proof}

We now apply these considerations to the affinoid $Y$ introduced in
\cite{mccallum:1988}, which is defined as follows.  Let $s$ and $t$ be
the functions on $F_{a,b,c}$ defined by
\begin{eqnarray}
  \label{eq:13}
  x &=& -\frac{a}{c}(1 + \pi^{(p-1)/2}s)\\
y &=& (-1)^ca^ab^bc^c(1+\pi t).
\end{eqnarray}
Let $Y$ be the affinoid defined over $L_w$ by the inequalities
\begin{equation*}
  |t| \leq |\pi_L^{-1}|,\quad |s| \leq 1.
\end{equation*}
A basis of holomorphic differentials on $F_{a,b,c}$ is
\begin{equation*}
  \omega_k = 
E_k \frac{x^{\left[\frac{ka}{p}\right]}(1-x)^{\left[\frac{kb}{p}\right]}}{y^k}
  \, dx, \quad k \in H_{a,b,c},
\end{equation*}
for some constants $E_k$ and where $H_{a,b,c}$ is a certain subset of
$\{1, 2, \dots , p-1\}$ of cardinality $(p-1)/2$ ($H_{a,b,c}$ can be
identified with the CM-type of $F_{a,b,c}$). We can and do choose the
constants $E_k$ so that $\omega_k$ has expansion
\begin{equation}
  \label{eq:14}
  \omega_k \equiv ds \quad \affmod{Y}{\pi_L},
\end{equation}
(note that this normalization is different from that of
\cite{mccallum:1988}).  Now, $P_1$ is the $\mclambda$-torsion point
represented by the divisor $(0,0) - \infty$, and we choose $f_1 = x$.
In \cite{mccallum:1988} it was shown that
\begin{equation}
  \label{eq:15}
\frac{df_1}{f_1} \equiv  \pi^{(p-1)/2} \sum_{k \in H_{a,b,c}} b_k
\omega_k \quad \affmod{Y}{p}
\end{equation}
for some $p$-adic integers $b_k$, satisfying
\begin{equation}
  \label{eq:16}
  \sum_{k \in H_{a,b,c}} b_k k^i \equiv
   \begin{cases}
     F & i=0\\
     0 & 1 \le i \le (p-3)/2
  \end{cases} \pmod{\pi_K}, \quad F \in \mathbb{Z}/p\mathbb{Z}^\times. 
\end{equation}
Note that although it was assumed that $F_{a,b,c}$ is wild split in
Section~5 of \cite{mccallum:1988}, there is nothing in the definition
of $Y$ or the calculation showing~(\ref{eq:15}) and~(\ref{eq:16}) that
uses this assumption. It is only at the end of that section that the
assumption comes in.

\begin{lemma}
  \label{thm:congruencelemma}
  If $$
  \sum_{k \in H_{a,b,c}} u_k \omega_k \equiv \sum_{k \in
    H_{a,b,c}} v_k \omega_k \quad \affmod{Y}{\pi^{n+ (p-3)/2}}$$
  then
  $$u_k \equiv v_k\pmod{\pi^n}, \quad k \in H_{a,b,c}.$$
\end{lemma}

\begin{proof}
  See pages 658--659 of \cite{mccallum:1988}.
\end{proof}

\begin{proposition}\label{thm:approximation}
  We have
  \begin{equation*}
\frac{df_{3}}{f_{3}} \equiv \sum_{k \in H_{a,b,c}}c_k \omega_k 
\affmod{Y}{p}, \qquad \text{$c_k \equiv 0 \pmod{\pi^{(p-5)/2}}$}
  \end{equation*}
  and
\begin{equation}
  \label{eq:17}
  \frac{df_{4}}{f_{4}} \equiv \sum_{k\in H_{a,b,c}} d_k
  \omega_k \affmod{Y}{p}, \qquad \text{$d_k \equiv
  -\pi^{(p-7)/2}\frac{b_k}{k^3} \pmod{\pi^{(p-5)/2}}$},
\end{equation}
where the $b_k$ are as in equation~\eqref{eq:15}.
\end{proposition}

\begin{proof}
  We have
  $$
  \mclambda^2_* \frac{df_{3}}{f_{3}} \equiv \frac {df_1}{f_1}
  \affmod{Y}{p}.
  $$
  Since $\zeta_* \omega_k = \zeta^{-k}\omega_k$, we have $\lambda_*
  \omega_k = \lambda^{\sigma} \omega_k$, for some $\sigma \in
  \Gal(K/\mathbb{Q}_p)$. Hence it follows from
  Lemma~\ref{thm:congruencelemma} and \eqref{eq:15} that
  $$\lambda^{2\sigma} c_k \equiv \pi^{(p-1)/2}b_k
  \pmod{\pi^{(p+1)/2}}.$$
  Thus $c_k \equiv 0 \pmod{\pi^{(p-5)/2}}$, as
  claimed.  Similarly, we have $\mclambda_*^3 (df_{4}/f_{4}) \equiv
  (df_1/f_1)$, so $\lambda^{3\sigma} d_k \equiv \pi^{(p-1)/2}b_k
  \pmod{\pi^{(p+1)/2}}$. Furthermore, since $\zeta^\sigma =
  \zeta^{-k}$, it follows from our choice of $\lambda$ that
  $\lambda^\sigma/\pi \equiv -k \pmod{\pi}$ for $1 \le k \le p-1$, so
  we get equation~\eqref{eq:17}.
\end{proof}

\begin{lemma}
  \label{thm:calculation}
  $$
  -\sum_{k \in H_{a,b,c}} \frac{b_k}{k^3} \equiv F( q(a^ab^bc^c)^3+
  a b c B_{p-3})\pmod{\pi}),
  $$
  where $F$ is as in~\eqref{eq:16}.
\end{lemma}

\begin{proof}
  Let $n = (p-1)/2$.  Define
  $$
  \Gamma_k(x_1, \dots, x_n) = \det \left[ \begin{array}{cccc} 1 & 1 & \dots & 1 \\
      x_1 & x_2 & \dots & x_n \\
      x_1^2 & x_2^2 & \dots & x_n^2 \\
      \vdots & \vdots & \ddots & \vdots \\
      x_1^{n-2} & x_2^{n-2} & \dots & x_n^{n-2}\\
      x_1^{n-1+k} & x_2^{n-1+k} & \dots & x_n^{n-1+k}
\end{array}
\right].
$$
Then an elementary
linear algebra calculation using~\eqref{eq:16} gives
\begin{equation*}
  \sum_{k \in H_{a,b,c}} \frac{b_k}{k^3} \equiv 
F \Gamma_3(H_{a,b,c}^{-1})/\Gamma_0(H_{a,b,c}^{-1})
\pmod{\pi}.
\end{equation*}
Let $S_i(x_1, \dots, x_n)$ be the $i$-th symmetric function. Then
$$
\Gamma_3 = \Gamma_0(S_1^3 - 2 S_1 S_2 + S_3).
$$
This can be proved by the usual method: the determinant vanishes if
$x_i = x_j$ for $i\neq j$, or if there is a polynomial of degree $n+2$
vanishing on the $x_i$, and with no term of degree $n-1$, $n$, or
$n+1$. Thus, if the roots of the polynomial are $x_1, \dots , x_n,
\alpha, \beta$, then
\begin{eqnarray*}
  \alpha + \beta + S_1 &=&0\\
S_2 + (\alpha+\beta) S_1 + \alpha \beta &=&0\\
\alpha \beta S_1 + \alpha S_2 + \beta S_2 + S_3 &=& 0
\end{eqnarray*}
Eliminating $\alpha$ and $\beta$ gives the condition $S_1^3 - 2S_1S_2
+ S_3 = 0$. Now, we have
\begin{eqnarray}
  S_1(H_{a,b,c}^{-1})  &\equiv & -q(a^ab^bc^c)\label{eq:18}\\
  S_2(H_{a,b,c}^{-1}) &\equiv& 0 \label{eq:19}\\
  S_3(H_{a,b,c}^{-1}) & \equiv & -\frac{B_{p-3}}{3}(a^3+b^3+c^3) \equiv -a b c B_{p-3}.\label{eq:20}
\end{eqnarray}
It is explained in \cite{mccallum:1988}, Lemma 5.24, how
equation~\eqref{eq:18} follows from \cite[17]{vandiver:1919};
equation~\eqref{eq:19} follows from parity considerations; and
equation~\eqref{eq:20} follows from \cite[16]{vandiver:1919},
in exactly the same way as \eqref{eq:18} follows from
\cite[17]{vandiver:1919}.
\end{proof}

We now define $p$-adic discs in $Y$, to which we apply
Lemma~\ref{thm:YtoX}. Let $X$ be the sub-affinoid of $Y$ defined by
$|t| \leq 1$ in the wild case and by $|s| \leq |\pi_K|$ in« the tame
case.  Let $E_w$ be the quadratic unramified extension of $K_w$.  If
$F_{a,b,c}$ is wild, $X$ is isomorphic to a union of two closed discs,
which are defined over $K_w$ in split case and over $E_w$ in the
non-split case.  Furthermore, $T = t$ is a parameter on each disc. If
$F_{a,b,c}$ is tame, then $X$ is isomorphic to a union of $p$ closed
discs defined over $K_w$, and $T = s/\pi_K$ is a parameter on each
disk. For proofs of these facts we refer the reader
to\cite{mccallum:1988} (where $T = s'$ in the tame case).  We denote
by $Z$ be any of the discs that are components of $X$, with parameter
$T$. We can write
\begin{equation*}
  f_i|_X = C_i u_i(T) v_i(T^{p}) g_i(T)^p, \quad i = 1, 2, 3, 4,
\end{equation*}
where $u_i$ and $v_i$ are unit power series with constant term 1 and
integer coefficients, $u_i$ has no terms in $T^{p}$, and $g_i$ is a
monic polynomial with integer coefficients. Furthermore, these
conditions uniquely determine the $u_i$, $v_i$ and $g_i$. Then
\begin{equation}
  \label{eq:21}
  \frac{df_i}{f_i} \equiv \frac{du_i}{u_i} \affmod{Z}{p}.
\end{equation}
For a $p$-adic field $H$ we denote by $U_H[[T]]$ the power series in
$\mathcal{O}_H[[T]]$ which are congruent to 1 modulo the maximal ideal
in $\mathcal{O}_H[[T]]$.

\begin{theorem}\label{thm:proposition-let-d}
  Let $Z$ be any of the discs that are components of $X$ and let $T$
  be a parameter on $Z$. Then
  \begin{equation}
    \label{eq:22}
    u_i = 1 + \pi^{(p+3)/2-i}D_i T + O(\pi^{(p+5)/2-i}T),
  \end{equation}
  where $|D_i| \leq 1$, $i = 1, 2, 3, 4$. Furthermore, $|D_1| = 1$,
  and under the hypotheses of Theorem~\ref{thm:nontrivial}, $|D_4| =
  1$ and
\begin{equation*}
  \frac{D_4}{D_1} \equiv 
    q(a^ab^bc^c)^3 + abc B_{p-3}.
\end{equation*}
Finally, $u_i$ for $i = 1, 2, 3$ are defined over $E_w$, and
\begin{equation*}
  u_4 \in 1 + \pi_K^{(p-5)/2}D_4 T U_{E_w}[[T]] + \pi_K^{(p+1)/2}\pi_L^{-3}ED_1
  T U_{F_w}[[T]],
\end{equation*}
where $E \in \mathbb{Z}/p\mathbb{Z}^\times$ is independent of the triple
$(a,b,c)$.
\end{theorem}

\begin{proof}
  In both wild and tame cases we have $\omega_k \equiv \pi D \, dT$
  $\affmod{Z}{\pi^2}$ for all $k \in H_{a,b,c}$, with $D \in
  \mathbb{Z}/p\mathbb{Z}^\times$ independent of $k$. This
  follows from our normalization~\eqref{eq:14}, since in the tame case
  we have
  \begin{equation}
 s = \pi T,\label{eq:23}
 \end{equation}
 and in the wild case it follows from \cite[(5.6)]{mccallum:1988},
 where it is shown that the expansion of $s$ in terms of $t$ on either
 of the discs in $X$ is
\begin{equation}
  \label{eq:24}
  s^2 = \frac{-q(a^ab^bc^c)2b}{ac} + \pi \frac{2b}{ac}(t^p - t) + O(\pi^2).
\end{equation} 
The statements about $D_1$ follow from~\eqref{eq:13}, \eqref{eq:23} (in
the tame case) and \eqref{eq:24} (in the wild case), since $f_1 = x$.
The statement about $D_2$ was proved in
\cite[Theorem~5.13]{mccallum:1988}.  Although this theorem is stated
only for the wild-split case, the consequence~\eqref{eq:22} is easily
seen to hold in the other cases as well (the part of Theorem~5.13
specific to the wild-split case translates into the statement $|D_2| =
1$ in the current notation, and we do not need it here). The
statements about $D_3$ and $D_4$ follow from
Proposition~\ref{thm:approximation}, Lemma~\ref{thm:calculation},
and~\eqref{eq:21}. The statement about the ratio $D_4/D_1$ follows
from~\eqref{eq:15}, the case $i=0$ of~\eqref{eq:16}, \eqref{eq:17},
and Lemma~\ref{thm:calculation}, taking note of the
normalization~\eqref{eq:14} and the fact that $ds \equiv
\mathrm{unit}\times \pi dT \affmod{Y}{\pi^2}$. The statements about
the fields of definition follow from the fact that $f_i$ is defined
over $K$ for $i = 1, 2, 3$ and $f_4$ is defined over $L$, and that the
discs $Z$ are always defined over $E_w$. The final statement follows
from considerations of ramification theory. Locally, we have
$\eta_{p-3} = 1 + a \pi^{p-3}$ modulo $p$-th powers, so the (upper)
conductor of $L_w/K_w$ is 3. Now, it follows from the properties of
the $P_i$ that $u_4^{\sigma-1} \equiv u_1$ modulo
$\mathcal{O}_{F_w}[[T]]^{\times p}U_{F_w}[[T^p]]$, and, since $u_1 \in 1 +
\pi^{(p+1)/2}TD_1 U_{F_w}[[T]]$, this implies the final statement with
$E$ such that $(\sigma -1)\pi_L^{-3} \equiv E^{-1} \pmod{\pi_L}$.
\end{proof}

\section{Computation of the Cassels pairing}
Recall the local descent maps
\begin{equation*}
\delta_i = \iota_{P_1}\times \cdots \iota_{P_i}: J(K_w) \rightarrow  
(K_w^\times/K_w^{\times p})^i
\end{equation*}
described in Section~\ref{sec:formulas-pairings}.  We start by
observing a couple of properties that follow from the choice of $P_i$
made in Section~\ref{sec:formulas-pairings}. First, we have
\begin{equation}
  \label{eq:25}
  \iota_{P_i}\circ \mclambda = \iota_{P_{i-1}}, \quad i = 2, 3, 4.
\end{equation}
Second, for $i = 1, 2, 3$ we have, from eigenspace considerations,
\begin{equation}
\iota_{P_i}(J(K_w)(k)) \subset V(k-i+1). \label{eq:26}
\end{equation}
Let $A \subset J(K_w)$ be the subgroup generated by divisors supported
on the discs $|T| \leq |\pi_K|$ in $X$.  Let
$$
V[i,j] = \bigcup_{i \le k \le j} V(i).
$$
Note that $V(i) = 0$ for $i>p$.
\begin{proposition}\label{thm:proposition-we-have}
  We have
\begin{equation}
 \iota_{P_i}(A)
  \subset V[(p+5)/2-i,p],\quad i = 1, 2, 3.\label{eq:27}
  \end{equation}
  If $F_{a,b,c}$ is wild split, we have
\begin{equation}
  \label{eq:28}
  \iota_{P_i}(J(K_w)) \subset V[(p+1)/2-i,p],\quad i = 1, 2, 3.
\end{equation}
If $F_{a,b,c}$ is wild non-split or tame, we have
\begin{equation}
\label{eq:29}
  \iota_{P_i}(J(K_w)) \subset V[(p+3)/2-i,p]\quad i = 1, 2, 3.
\end{equation}
Furthermore, in the case $i=1$, the containments in~\eqref{eq:27}
and~\eqref{eq:29} are equalities.
\end{proposition}

\begin{proof}
  The containments in~\eqref{eq:27} follow immediately from
  Proposition~\ref{thm:proposition-let-d}, as does the claim that the
  inclusion is equality in the case $i=1$. Now
  Faddeev\cite{faddeev:1961} proved that
  \begin{equation*}
    \im \iota_{P_1} = 
    \begin{cases}
      V((p-1)/2) \cup V[(p+3)/2, p] & \text{$F_{a,b,c}$ is wild
        split}\\
      V[(p+3)/2,p] & \text{otherwise}
    \end{cases}
  \end{equation*}
  This implies the statements~\eqref{eq:28} and~\eqref{eq:29} in the
  case $i=1$, and also that the image of $A$ in $J(K_w)/\mclambda
  J(K_w)$ has codimension 1, and the eigenvalue of the quotient is
  $\kappa^{(p-1)/2}$ in the wild non-split case and $\kappa^{(p+3)/2}$
  in the other cases. The remaining statements now follow from
  \eqref{eq:25}, \eqref{eq:26}, and~\eqref{eq:27}.
\end{proof}

\begin{proposition}\label{thm:proposition-if-f_a-1}
  If $F_{a,b,c}$ is wild non-split or tame, then $$\delta_2(J(K_w)) =
  V[(p+1)/2,p]^2.$$
\end{proposition}

\begin{proof}
  From~\eqref{eq:25} with $i=2$ we have
  \begin{equation}
\im \delta_2 \cap (K^\times/K^{\times p}) \times 1  =  \im
  \delta_1 \times 1 = V[(p+1)/2,p] \times 1\label{eq:30}
  \end{equation}
  Furthermore, given $u \in V[(p+3)/2, p]$, we can find $a \in A$ such
  that $\iota_{P_1}(a) = u$, and $\iota_{P_2}(a) \in V[(p+1)/2,p]
  \subset \im \iota_{P_1}$.  Thus, modifying $a$ by $\mclambda
  J(K_w)$, we can choose it so that $\iota_{P_2}(a) = 1$. Hence
  \begin{equation}
    \label{eq:31}
\im  \delta_2 \supset 1 \times V[(p+3)/2,p]. 
  \end{equation}
  Now, it follows from local duality that $\im \delta_2$ must be
  maximal isotropic with respect to the cup product pairing on
  $(K^\times/K^{\times p})^2 = H^1(K, J[\mclambda^2])$ induced by the
  Weil pairing on $J[\mclambda^2]$. Since $\lambda^2 = \hat\lambda^2$,
  the Weil pairing is skew symmetric. Thus the pairing on
  $(K^\times/K^{\times p})^2$ is a non-zero multiple of $((a_1, b_1),
  (a_2, b_2)) \mapsto (a_1, b_2)_w(b_1, a_2)_w^{-1}$, where $(,)_w$
  denotes the Hilbert symbol at $w$. The only maximal isotropic
  subgroup satisfying~\eqref{eq:30} and~\eqref{eq:31} is the one given
  in the statement of the proposition.
\end{proof}
Define a subspace $V_{\mathrm{global}} \subset V$ by 
\begin{equation*}
  V_{\mathrm{global}} = \bigoplus_{\text{$i$ even, $2 \leq i \leq
  p-3$}} V(i).
\end{equation*}

\begin{proposition}\label{thm:proposition-assume-p}
  Assume $p \geq 11$ and $F_{a,b,c}$ is wild non-split or tame.  There
  exists a point $x \in A$ such that $\iota_{P_1}(x)$ generates $V((p+5)/2)$
  and $\iota_{P_i}(x) \in V_{\mathrm{global}}$ for $i = 2, 3$.
\end{proposition}

\begin{proof}
  It follows from Proposition~\ref{thm:proposition-we-have} that $A$,
  regarded as a $Z_p[\zeta]$-submodule of $J(K_w)$, has codimension at
  most $1$. Hence $\delta_3(A)$ has codimension at most 3 as a
  $\mathbb{F}_p$-vector space in $\delta_3(J(K_w))$. By
  Proposition~\ref{thm:proposition-we-have} we can choose $x \in A$
  such that $\iota_{P_1}(x)$ generates $V((p+5)/2)$. This condition 
 leaves freedom to modify $x$
  by anything in $\mclambda J(K_w)$, which would change $\delta_3(x)$
  by anything in $\im \delta_2$. Thus, modifying $x$ as needed, we can
  ensure that $\iota_{P_i}(x) \in V_{\mathrm{global}}$, $i = 2, 3$.
  The number of degrees of freedom in performing this modification is
  equal to the dimension of $\im \delta_2 \cap V_{\mathrm{global}}$,
  which is at least 4 if $p\geq 11$, by
  Proposition~\ref{thm:proposition-if-f_a-1}. Thus we can ensure that
  $x$ remains in $A$ when making the modification.
\end{proof}

\begin{proof}[Computation of Cassels pairing for Theorem~\ref{thm:free}]
  We show here that the pairing $\langle , \rangle_2$ is trivial under
  the hypotheses of Theorem~\ref{thm:free}. In the next section we
  explain how this implies the theorem. 

Denote by $\ell_i:
  S_{\mclambda^i} \to J(K_w)/\mclambda^i J(K_w)$ the localization map.
  We claim that, under the hypotheses of Theorem~\ref{thm:free},
  $\iota_{P_1}(\ell_1(S_\mclambda)) \subset V[(p+3)/2,p]$ or
  $\iota_{P_1}(\ell_1(S_\mclambda)) \subset V((p+1)/2) \cup
  V[(p+5)/2,p]$. Now, $V(i)$ pairs non trivially with $V(j)$ under
  Hilbert pairing if and only if $i + j \equiv p \pmod{p-1}$. Thus, it
  follows from our claim and from~\eqref{eq:26} that
  $\iota_{P_3}(\ell_3(J(K_w))$ pairs trivially with
  $\iota_{P_1}(\ell_1(J(K_w))$.
  
  To see the claim, note that if hypothesis (\ref{item:1}) of
  Theorem~\ref{thm:free} is satisfied, namely that $F_{a,b,c}$ is
  wild-split and $p \equiv 3 \pmod{4}$, then, by~\eqref{eq:29},
  $\ell_1(S_\mclambda) \subset V((p-1)/2)\cup V[(p+3)/2,p]$.
  Furthermore, we can eliminate $V((p-1)/2)$ as a possibility, because
  $\ell_w$ factors through $H^1(K(p)/K, \mu_p) \rightarrow H^1(K_w,
  \mu_p)$. Since $(p-1)/2$ is odd, it follows from~\eqref{eq:7} that
  $\ell_w$ can have nontrivial image in $V((p-1)/2)$ only if
  $C((p-1)/2)$ is nontrivial, which would imply $p \mid B_{(p+1)/2}$.
  This never happens if $p \equiv 3 \pmod{4}$.
  
  If hypothesis (\ref{item:2}) of Theorem~\ref{thm:free} is
  satisfied, namely that $F_{a,b,c}$ is wild non-split or tame and the
  image of $U$ in either $V((p+1)/2)$ or $V((p+3)/2)$ is trivial, then
  the claim follows immediately from~\eqref{eq:29}.
\end{proof}

\begin{proof}[Proof of Theorem~\ref{thm:nontrivial}]
  We exhibit $a \in S_{\lambda^3}$ and $b \in S_\lambda$ which pair
  nontrivially under the Cassels pairing. 

Recall that $S_{\lambda^i}
  \subset H^1(K(p)/K, J[\lambda^i])$, and since $p$ is regular, this
  latter group is isomorphic to
  $(\mathcal{O}_K^\times/\mathcal{O}_K^{\times p})^i$ for $i \leq 3$.
  The Selmer group is the subgroup obtained by imposing the local
  conditions at $w$. Since $(p+1)/2$ is even, we can choose an element
  $b \in \mathcal{O}_K^\times/\mathcal{O}_K^{\times p}$ which
  generates $V((p+1)/2)$, and $b$ satisfies
  the local condition by Proposition~\ref{thm:proposition-we-have}, so
  $b \in S_\lambda$. 
  
  As for $a$, by Proposition~\ref{thm:proposition-assume-p} there
  exists $a_w = (a_{w,1}, a_{w,2}, a_{w,3}) = \delta_3(x)$, $x \in A$,
  such that $a_{w,1}$ generates $V((p+5)/2)$ and $a_{w,2}, a_{w,3} \in
  V_{\mathrm{global}}$.  Choose $a_i \in
  \mathcal{O}_K^\times/\mathcal{O}_K^{\times p}$ specializing to
  $a_{w,i}$ for $i = 1, 2, 3$ and define $a \in S_{\lambda^3}$ by $a =
  (a_1, a_2, a_3)$. 
  
  Now, by Lemma~\ref{thm:lemma-suppose-that}, $\mclambda_*^2 a = a_1
  \in V((p+5)/2)$ is the norm of a global unit $\epsilon$ in
  $\mathcal{O}_L^\times$, and by Proposition~\ref{thm:prop-supp-p},
s the Cassels pairing of $a$ and $b$ is the Hilbert pairing
  $(c_w, b_w)$, where
\begin{equation}
  \label{eq:32}
  c_w = \iota_{P_4}(x)^{-1}  N'\epsilon
  \quad \text{in $L_w^\times/L_w^{\times p}$}.
\end{equation}
To prove that the pairing is nontrivial, it suffices to show that
$c_w$ is not a $p$-th power, and for that it suffices to show that its
image in $L_w^\times/L_w^{\times p}$ is nontrivial. We may assume
without loss of generality that $c_w$, $N'\epsilon$, and
$\iota_{P_4}(x)$ are eigenvectors for a lift $\tilde \Delta$ of
$\Delta = \Gal(K_w/\mathbb{Q}_p)$ to $\Gal(L_w/\mathbb{Q}_p)$, with
eigenvalue $\kappa^{(p-1)/2}$.

Since $x \in A$, we may choose a divisor $D$ supported on $|T|\leq
|\pi|$ such that $a_{w,i} = f_i(D)$, $1 \leq i \leq 3$. Since $D$ is
supported on $|T| \leq |\pi|$, we have
  \begin{equation*}
u = f_4(D) \equiv u_4(D) \pmod{L_w^{\times p}(1 +
  \pi^p\mathcal{O}_{L_w})}.
\end{equation*}
From the galois properties of the $P_i$, we have
\begin{equation}\label{eq:33}
u \in (L_w^\times/L_w^{\times
  p})((p-1)/2)), \quad (\sigma-1)u = v,
\end{equation}
where $v$ is the image in $L_w^\times/L_w^{\times p}$ of a generator
of $V((p+5)/2)$.  Since $p \geq 19$, $(p+5)/2$ is less than $p-3$, and
thus $v \neq 0$. Thus the subspace of $L_w^\times/L_w^{\times p}$
satisfying the conditions~\eqref{eq:33} is two-dimensional, with
generators $u_1$ and $u_2$, where $u_1$ is the image of a generator of
$V((p-1)/2)$ with expansion $u_1 = 1 + \pi_K^{(p-1)/2} +
O(\pi_K^{(p+1)/2})$, and $u_2 \in (L_w^\times/L_w^{\times
  p})((p-1)/2)$ has expansion $u_2 = 1 + \pi_K^{(p+5)/2}\pi_L^{-3} +
O(\pi_K^{(p+5)/2})$.  Thus $u = \alpha u_1 + \beta u_2$ for some
$\alpha, \beta \in \mathbb{Z}/p\mathbb{Z}$.  Comparing with
Proposition~\ref{thm:proposition-we-have} we see that
  \begin{equation}
    \label{eq:34}
\alpha/\beta = \gamma(a,b,c) =
    q(a^ab^bc^c)^3 + abc B_{p-3}. 
  \end{equation}
  Now, we may replace $(a,b,c)$ by any $(a', b', c') \equiv (ta, tb,
  tc) \pmod{p}$, for $t \in \mathbb{F}_p^\times$. It is easily seen,
  using the property $q(xy) \equiv q(x) + q(y)$, that $\gamma(ta, tb,
  tc) \equiv t^3 \gamma(a,b,c)$. Thus, from~\eqref{eq:34}, we see that
  by varying $t$ appropriately we may ensure that $u$, and hence
  $c_w$, varies in $L_w^\times/L_w^{\times p}$, and, in particular,
  takes on nonzero values. Hence there exists a choice of $t$ such
  that the pairing is nontrivial for the curve $F_{a', b', c'}$.
  However, this curve is isomorphic to $F_{a,b,c}$, and hence the
  pairing must be nontrivial in that case as well.
\end{proof}

\section{Shafarevich-Tate groups}\label{sec:shaf-tate-groups}
The proofs of Theorem~\ref{thm:free} and
Theorem~\ref{prop:exact-p-primary} follow from the computations of
the Cassels-Tate pairing by means of the following proposition.
\begin{proposition}\label{thm:prop-all-posit}
  For all positive integers $m$ and $n$, the restriction of the
  Cassels-Tate pairing induces a perfect pairing
  $$(\Sh[\la^m]/(\la^n\Sh[\la^{n+m}])) \ \times (\Sh[\la^n]/(\la^m
  \Sh[\la^{n+m}])) \ \longrightarrow \ \Q/\Z.$$
\end{proposition}
\medskip

Let $\Sh_{\div}$ denote the maximal divisible subgroup of $\Sh$, i.e.
$x \in \Sh_{\div}$ if and only if for every non-zero integer $n$ there
exists $y \in \Sh$ such that $x = n y$. Let $\Sh_{\red}$ denote the
quotient group $\Sh/\Sh_{\div}$. Note that:

\begin{lemma}\label{thm:lemma-sh_d-divis}
  $\Sh_{\div}$ is a divisible group in the usual sense, i.e.\ for
  every non-zero $n \in \Z$ multiplication by $n$ on $\Sh_{\div}$ is
  surjective.
\end{lemma}

\begin{proof}
  The argument is standard: Since $\Sh[m]$ is finite for all non-zero
  $m \in \Z$, the groups $N\Sh[Nm]$, $N>0$, stabilize for sufficiently
  large $N$. Thus for every $m$ there is an integer $N(m)$ such that
  if an element of $\Sh[m]$ is divisible by $N(m)$ it is infinitely
  divisible. Now if $x \in \Sh_{\div}[m]$ and $n>0$, choose $y \in
  \Sh[N(nm)nm]$ such that $N(nm)ny = x$. Then $y' = N(nm)y$ is in
  $\Sh_{\div}[nm]$ and $ny' = x$.
\end{proof}

Note that since $\zeta$ is an automorphism of $\Sh$ it preserves
$\Sh_{\div}$, and hence so does $\Z[\zeta]$. Furthermore, since
$\mclambda^{p-1}$ is a unit times $p$ in $\Z[\zeta]$, $\Sh_{\div}$ is
divisible by $\mclambda^{n}$ for any positive $n$.

\begin{lemma}\label{thm:lemma-exact-sequence}
  The exact sequence
  $$0 \longrightarrow \Sh_{\div} \longrightarrow \Sh \longrightarrow
  \Sh_{\red} \longrightarrow 0 $$
  induces by restriction an exact
  sequence
  $$0 \longrightarrow \Sh_{\div}[\la^n] \longrightarrow \Sh[\la^n]
  \longrightarrow \Sh_{\red}[\la^n] \longrightarrow 0 $$
  for any
  positive integer $n$.
\end{lemma}

\begin{proof}
  Only the surjectivity is in question.  Let $x \in
  \Sh_{\red}[\la^n]$. Lift $x$ to $y \in \Sh$.  Then $\la^n y = z \in
  \Sh_{\div}$. By Lemma~\ref{thm:lemma-sh_d-divis}, we can find $w \in
  \Sh_{\div}$ such that $\la^n w = z = \la^n y$. But then $y-w \in
  \Sh[\la^n]$ and $y-w$ reduces to $x$ in $\Sh_{\red}$.
\end{proof}

It is well-known that $\Sh_{\red}[p^{\infty}]$ is a finite group and
that the Cassels-Tate pairing induces a perfect pairing
$$[\cdot \ , \cdot]: \Sh_{\red}[p^{\infty}] \ \times \ 
\Sh_{\red}[p^{\infty}] \ \longrightarrow \ \Q/\Z. $$
We now have the
following lemma:

\begin{lemma}\label{thm:lemma-annih-sh_r}
  The annihilator of $\Sh_{\red}[\la^m]$ with respect to the latter
  pairing equals $\la^m \Sh_{\red}[p^{\infty}]$, for all positive
  integers $m$.
\end{lemma}

\begin{proof}
  It is clear from the definition of the pairing given in
  \cite{mccallum:1988}, for example, and from the functoriality
  properties of the Weil pairing, that $[\zeta a , a'] = [a,
  \zeta^{-1} a']$. Hence, if $\hat\la = \zeta^{-1} -1$, then
  $\hat\la^m \Sh_{\red}[p^{\infty}]$ annihilates $\Sh_{\red}[\la^m]$.
  Since $\hat \mclambda/\mclambda$ is a unit in $\Z[\zeta]$, we have
  $\hat\la^m \Sh_{\red}[p^{\infty}] = \la^m \Sh_{\red}[p^{\infty}]$.
  So the kernel $H$ on the right factor of the restricted pairing
  $$\Sh_{\red}[\la^m] \times \Sh_{\red}[p^{\infty}] \longrightarrow
  \Q/\Z $$
  contains $\la^m \Sh_{\red}[p^{\infty}]$. Note that the
  kernel on the left factor of the latter pairing is trivial.
  Therefore, the cardinalities of $\Sh_{\red}[\la^m]$ and
  $\Sh_{\red}[p^{\infty}]/H$ are equal. But
  $$|\Sh_{\red}[p^{\infty}]| = |\Sh_{\red}[\la^m]| \ |\la^m
  \Sh_{\red}[p^{\infty}]|, $$
  hence $H=\la^m \Sh_{\red}[p^{\infty}]$.
\end{proof}

\begin{lemma}\label{thm:lemma-all-positive}
  For all positive integers $m$ and $n$, the restriction of the
  Cassels-Tate pairing induces a perfect pairing
  $$(\Sh_{\red}[\la^m]/(\la^n \Sh_{\red}[\la^{n+m}])) \ \times \ 
  (\Sh_{\red}[\la^n]/( \la^m \Sh_{\red}[\la^{n+m}])) \ \longrightarrow
  \ \Q/\Z. $$
\end{lemma}

\begin{proof}
  By Lemma~\ref{thm:lemma-annih-sh_r}, the annihilator of
  $\Sh_{\red}[\la^m]$ in $\Sh_{\red}[\la^n]$ equals
  $$\la^m \Sh_{\red}[p^{\infty}] \ \cap \ \Sh_{\red}[\la^n] = \la^m
  \Sh_{\red}[\la^{n+m}], $$
  and the assertion follows.
\end{proof}

\begin{proof}[Proof of Proposition~\ref{thm:prop-all-posit}.]
  By Lemma~\ref{thm:lemma-all-positive}, it suffices to show that for
  all $m$ and $n$ the groups $\Sh[\la^m]/(\la^n \Sh[\la^{n+m}])$ and
  $\Sh_{\red}[\la^m]/(\la^n \Sh_{\red}[\la^{n+m}])$ are isomorphic.
  By Lemma~\ref{thm:lemma-exact-sequence}, we have a commutative
  diagram
  $$
  \begin{array}{lllllllll} 0 & \longrightarrow &
    \Sh_{\div}[\la^{n+m}] & \longrightarrow & \Sh[\la^{n+m}] &
    \longrightarrow & \Sh_{\red}[\la^{n+m}] &
    \longrightarrow & 0 \\
    \\
    & & \downarrow \alpha =\la^n & & \downarrow \beta=\la^n & &
    \downarrow  \gamma=\la^n \\
    \\
    0 & \longrightarrow & \Sh_{\div}[\la^m] & \longrightarrow &
    \Sh[\la^m] & \longrightarrow & \Sh_{\red}[\la^m] & \longrightarrow
    & 0
\end{array} $$
where the horizontal sequences are exact.  By the snake lemma, we get
an exact sequence
$$
0 \to Ker(\alpha) \to Ker(\beta) \to Ker(\gamma) \to Coker(\alpha)
\to Coker(\beta) \to Coker(\gamma) \to 0.
$$
By Lemma~\ref{thm:lemma-sh_d-divis}, we have $Coker(\alpha) =0$,
hence $Coker(\gamma)$ is isomorphic to $Coker(\beta)$, and this
completes the proof.
\end{proof}

\begin{proof}[Proof of Theorem~\ref{thm:free}.]
  By the structure theorem for torsion modules over Dedekind domains
  we have a $\Z[\zeta]$-module decomposition
  $$\Sh[\la^3] \iso (\Z[\zeta]/(\la))^{t_1} \ \oplus \ 
  (\Z[\zeta]/(\la^2))^{t_2} \ \oplus \ (\Z[\zeta]/(\la^3))^{t_3} , $$
  where $t_1$, $t_2$ and $t_3$ are non-negative integers.  The
  computations in the previous section show that the pairing (obtained
  by restricting the Cassels-Tate pairing)
  $$\Sh[\la^2] \ \times \ \Sh[\la] \ \to \ \Q/\Z $$
  is trivial.  By
  Proposition~\ref{thm:prop-all-posit} (for $m=2$ and $n=1$), we get
  that the groups $\Sh[\la^2]/(\la\Sh[\la^3])$ and
  $\Sh[\la]/(\la^2\Sh[\la^3])$ are both trivial. But then
  $$(\Z[\zeta]/(\la))^{t_1} \oplus \ (\Z[\zeta]/(\la))^{t_2} \oplus
  (\Z[\zeta]/(\la))^{t_3} \iso \Sh[\la]=\la^2 \Sh[\la^3] \iso
  (\Z[\zeta]/\la)^{t_3} $$
  so $t_1=t_2=0$, which proves the claim.
\end{proof}

\begin{proof}[Proof of Theorem~\ref{prop:exact-p-primary}.] 
  Let
  $$\Sh[\la^4] \iso (\Z[\zeta]/(\la))^a \ \oplus \ 
  (\Z[\zeta]/(\la^2))^b \ \oplus \ (\Z[\zeta]/(\la^3))^c \ \oplus \ 
  (\Z[\zeta]/(\la^4))^d. $$
  If we show that $d=0$, then $\la^3$
  annihilates $\Sh[\la^4]$, therefore $\Sh[\la^4]=\Sh[\la^3]$. By
  induction, this implies
  $\Sh[p^{\infty}]=\Sh[\la^{\infty}]=\Sh[\la^3]$. So assume $d \ge 1$.
  Since the Cassels-Tate pairing on $\Sh[\la^3] \times \Sh[\la]$ is
  non-trivial, Proposition~\ref{thm:prop-all-posit} implies that
  $\Sh[\la^3]/(\la \Sh[\la^4])$ has dimension $\ge 2$ over $\F_p$. Now
  $$\la \Sh[\la^4] \iso (\Z[\zeta]/(\la))^b \ \oplus \ 
  (\Z[\zeta]/(\la^2))^c \ \oplus \ (\Z[\zeta]/(\la^3))^d. $$
  Counting
  $\F_p$-dimensions, we get $6-(b+2c+3d) \ge 2$, therefore $b+2c+3d
  \le 4$. This implies $d=1$ and $c=0$. Therefore,
  $$\Sh[\la^4] \iso (\Z[\zeta]/(\la))^a \ \oplus \ 
  (\Z[\zeta]/(\la^2))^b \ \oplus \ (\Z[\zeta]/(\la^4)). $$
  This
  implies that
  $$(\Z[\zeta]/(\la))^2=\la^2(\Z[\zeta]/(\la^3))^2 \iso
  \la^2\Sh[\la^3] \subseteq \la^2\Sh[\la^4] \iso \Z[\zeta]/(\la^2) ,
  $$
  a contradiction.
\end{proof}

\section{Tame reduction}\label{sec:tame-reduction}

Although it is not strictly necessary for
Theorem~\ref{cor:diophant}, we take the opportunity to prove a
general lemma on tame reduction, since it clears up some confusion in
the literature. In \cite{lim:1995}, an attempt was made to improve the
result of \cite{mccallum:1988} on the existence of non-trivial
elements in $\Sh[\la]$ in the wild-split case, under the additional
hypothesis that the Jacobian of the Fermat curve in question is
non-simple.  However, as Lemma~\ref{thm:lemma-let-a} shows,
non-simple Jacobian and wild split reduction over $\Z_p[\zeta]$
are incompatible properties, so  the
Mordell-Weil rank estimates given in the last section of
\cite{lim:1995} are incorrect. As far as we can tell, the problem lies
in the use of the function $q(x)$ which computes the reduction type
(see the introduction).  Here as well as in \cite{mccallum:1988}, $q$
is evaluated on triples $(a,b,c)$ of integers such that $0 < a, b ,
a+b< p$ and $a+b+c=0$.  In \cite{lim:1995} however, $q$ is evaluated
on triples $(a,b,c)$ such that $0 < a, b, a+b < p$ and $a+b+c=p$.
While it does not make any difference which of the two types of
triples one chooses to define the curve $F_{a,b,c}$, it does make a
difference which type of triple one uses to evaluate $q$ and hence the
reduction type.  We have the following lemma:

\begin{lemma}\label{thm:lemma-let-a}
  Let $(a,b,c)$ be such that $J_{a,b,c}$ is non-simple.  Then
  $F_{a,b,c}$ has tame reduction over $\Z_p[\zeta]$.
\end{lemma}

\begin{proof}
  By \cite{koblitz-rohrlich:1978}, $J_{a,b,c}$ is non-simple if and
  only if $p \equiv 1 \pmod{3}$ and $F_{a,b,c}$ is isomorphic to
  $F_{1,r,-(r+1)}$, where $r^2+r+1=0$ in $\F_p$.  By definition of
  $q(x)$, it therefore suffices to show that $(r+1)^{(r+1)(p-1)} -
  r^{r(p-1)} \equiv 0 \pmod{p^{2}}$.  Since 6 divides $p-1$, it
  suffices to show that
  $$(r+1)^{6(r+1)}- r^{6r} \equiv 0 \pmod{p^2}. $$
  Let $k$ be an
  integer such that $r^2+r+1=pk$. Then $(r+1)^2=pk+r$. Therefore,
  $$(r+1)^6=(pk+r)^3 \equiv (r^3+3 r^2 pk) \pmod{p^2}. $$
  Hence,
  $(r+1)^{6(r+1)} \equiv (r^3+3r^2pk)^{r+1} \equiv (r^{3(r+1)} + 3r^2
  pk(r+1)r^{3r}) \pmod{p^2}$.  Now note that, since $r$ is a cube root
  of unity modulo $p$, $r^{3r}r^2 (r+1) \equiv -r \pmod{p}$, so $3r^2
  pk (r+1)r^{3r} \equiv -3rpk \pmod{p^2}$.  Hence, $(r+1)^{6(r+1)}
  \equiv (r^{3r+3} - 3rpk) \pmod{p^2}$.  Therefore,
  $$(r+1)^{6(r+1)} - r^{6r} \equiv (r^{3r}(r^3-r^{3r})-3rpk)
  \pmod{p^2}. $$
  Since $r^3=pk(r-1)+1$, we get $r^{3r} \equiv
  (rpk(r-1)+1) \pmod{p^2}$, so $r^3-r^{3r} \equiv -pk(r-1)^2
  \pmod{p^2}$. Hence,
  $$(r+1)^{6(r+1)} - r^{6r} \equiv -pk (r^{3r}(r-1)^2 +3r) \pmod{p^2}.
  $$
  Since $r^{3r}(r-1)^2+3r \equiv 0 \pmod{p}$, this proves the
  proposition.
\end{proof}

\begin{proof}[Proof of Theorem~\ref{cor:diophant}] 
  By Lemma~\ref{thm:lemma-let-a}, the reduction is tame in this case.
  By Theorem~\ref{thm:nontrivial} and
  Proposition~\ref{thm:prop-all-posit}, the $\F_p$-dimension of
  $\Sh[\la]/(\la^3\Sh[\la^4])$ is $\ge 2$. In particular, the
  $\F_p$-dimension of $\Sh[\la]$ is $\ge 2$.  Since $p$ is regular,
  the results of Faddeev (\cite{faddeev:1961}) show that the Selmer
  group $S_{\la}$ is 3-dimensional over $\F_p$. On the other hand,
  Gross and Rohrlich (\cite{gross-rohrlich:1978}) have shown that the
  Mordell-Weil rank of $J_{7,1,-8}$ over $\Q$ is non-zero. Therefore,
  the rank equals 1 and $\Sh[\la]$ is 2-dimensional over $\F_p$. By
  Theorem~\ref{thm:free} it follows that $\Sh[\la^3]$ has rank 2 over
  $\Z[\zeta]/(\la^3)$. Theorem~\ref{prop:exact-p-primary} then
  implies that $\Sh[p^{\infty}] \iso (\Z[\zeta]/(\la^3))^2$.  The
  statement about quadratic points on $F_{7,1,-8}$ and on the Fermat
  curve $X^{19}+Y^{19}+Z^{19}=0$ follows immediately from Corollary
  2.2 and Theorem 1.3 of \cite{tzermias:2002}.
\end{proof}

\bibliographystyle{amsalpha} 

\providecommand{\bysame}{\leavevmode\hbox to3em{\hrulefill}\thinspace}
\providecommand{\MR}{\relax\ifhmode\unskip\space\fi MR }
\providecommand{\MRhref}[2]{%
  \href{http://www.ams.org/mathscinet-getitem?mr=#1}{#2}
}
\providecommand{\href}[2]{#2}

\end{document}